\documentstyle[12pt,twoside]{article}
\newtheorem{theorem}{Theorem}

\newtheorem{lemma}[theorem]{Lemma}

\newtheorem{definition}[theorem]{Definition}

\begin{document}

%\begin{document}
\overfullrule=0pt
\baselineskip=24pt
% if i uncomment the preceeding line, the text will be double-spaced.
\font\tfont= cmbx10 scaled \magstep3
\font\sfont= cmbx10 scaled \magstep2
\font\afont= cmcsc10 scaled \magstep2
\title{\tfont Global Structures on CR Manifolds via Nash 
Blow-ups }
\bigskip
\author{ Thomas Garrity \\
\\ Department of Mathematics\\ Williams College\\Williamstown, MA  
01267\\ 
email:tgarrity@williams.edu\\
\bigskip \\
 Dedicated to William Fulton on his sixtieth 
birthday}

\date{}

\maketitle
\begin{abstract}
% Abstract text begins here

A generic compact real codimension two submanifold $X$ of  
$C^{n+2}$ will have a CR structure at all but a finite number of 
points (failing at the complex jump points $\cal {J}$).  
The main theorem of this paper gives a method of extending 
the CR structure on the non-jump points $X-\cal {J}$ to the jump points. We 
examine a Gauss map from $X-\cal {J}$ to an appropriate flag manifold $F$ and 
take the closure of the graph of this map in $X \times F$.  This is a 
version of a Nash blow-up.  We give a 
clean criterion for when this closure is a smooth manifold and see that 
the local differential properties at the points $X-\cal {J}$ can now be 
naturally extended to this new smooth manifold, allowing global 
techniques from differential geometry to be applied to compact 
CR manifolds.  As an example, we find topological obstructions for the 
manifold to be Levi nondegenerate.

\end{abstract}

\section{Introduction}
Let $X$ be a compact real $2n+2$ dimensional submanifold of the complex 
space $C^{n+2}$. For generic such $X$, at all but a finite number of 
points, the tangent space of $X$ will have a $2n$ dimensional subspace 
$H$
that inherits a complex structure  from the ambient 
$C^{n+2}$.  There are, though, topological obstructions preventing the 
subspaces $H$ from forming a subbundle of  the tangent bundle $TX$.  
The existence of such obstructions was shown by Wells \cite{Wells}. 
Lai \cite{Lai} gave an explicit description of these obstructions.

 There has recently been a lot of work on determining when two CR 
 structures are locally equivalent, subject to various restrictions 
 on dimension and conditions on the Levi form. There is the work of 
 Beloshapka \cite{Beloshapka1}, 
\cite{Beloshapka2}, Ebenfelt \cite{Ebenfelt}, \cite{Ebenfelt2}, Ezhov 
and Isaev \cite{Ezhov-Isaev}, Ezhov, Isaev and Schmalz 
\cite{Ezhov-Isaev-Schmalz}, Ezhov and Schmalz \cite{Ezhov-Schmalz1}, 
\cite{Ezhov-Schmalz2}, \cite{Ezhov-Schmalz3},
\cite{Ezhov-Schmalz4}, Garrity and Mizner \cite{GM}, \cite{GM2},
Le \cite{Le},  Mizner  \cite{Mizner} and  
Schmalz and Slovak \cite{Schmalz-Slovak}.  These works 
concentrate on the understanding of the Levi form, a vector-valued 
Hermitian form at each point mapping $H \times H$ to $TX/H$.

All of these techniques and methods for producing local invariants
 break down for compact manifolds.  
 What has  prevented people from applying standard 
tools from differential geometry to understand the obstructions 
preventing the extensions of these local invariants to global 
invariants has been that the subbundle $H$ is not a true subbundle.  
All of the local calculations depend on $H$, the part of the tangent 
bundle inheriting a complex structure from $C^{n+2}$, having real 
dimension $2n$.  For a compact $X$, there will be points (the {\it 
complex jump points}, which we will denote by $\cal {J}$) where the $H$ 
will have real dimension $2n+2$.   The existence of these points is 
what prevents any easy attempt to extend local invariants to 
global ones.

We use a version of the Nash blow-up to replace $X$, subject to 
certain natural conditions, with a smooth 
manifold $\tilde{X}$ so that there is a natural map $\pi: \tilde{X} 
\rightarrow X$ with $\pi$ an isomorphism from $\tilde{X} - 
\pi^{-1}\cal {J}$ 
to $X-\cal {J}$ and so that there 
is a complex rank $n$ vector bundle $\tilde{H}$ on 
$\tilde{X}$ such that $\tilde{H}$ pushes forward to the bundle $H$ on 
$X-\cal {J}$.  Thus global calculations can now be  performed.

The method presented here is to show that there is a natural map 
(a version of the 
Gauss map) from $X-\cal {J}$ to a flag manifold $F$.  The Nash blow-up is the 
closure of the graph of this map in $X \times F$.  Our main result is 
to give a clear criterion as to when this closure is a smooth 
manifold.  We will show that the Nash blow-up will be smooth when
 the Gauss-Lai image of $X$ transversally intersects the 
subvariety of real $2n$-planes in the real Grassmannian $G(2n, C^{n+2})$
that inherit a complex structure from the ambient $C^{n+2}$.

Finally, it gives me great pleasure to present this paper in honor of 
William Fulton's sixtieth birthday.

\section{Basic Definitions}

\subsection{CR structures}

Let $X$ be a compact real codimension two submanifold of $ C^{n+2}$.
Thus $X$ has real dimension $2n+ 2$.  Let $J:C^{n+2} \rightarrow 
C^{n+2} $ be the linear map corresponding to multiplication by $i$.  
Thus $J^{2} = -I$. For more on this, see \cite{Boggess}, chapter 
three, \cite {Chirka}, \cite {Jacobowitz}, \cite{Taiani} and 
\cite{Tumanov}.

\begin{definition}
The complex tangent space of X at a point p is the subspace
$$H_{p} = T_{p}X \bigcap JT_{p}.$$
\end{definition}
The complex tangent space is the subspace of the tangent space that 
inherits a complex structure from the ambient complex space 
$C^{n+2}$.  As we will discuss,  at all but
 a finite number of 
points for generic $X$, 
the real dimension of the complex tangent space $H_{p}$ will 
be $2n$ and thus complex dimension will be $n$.  

\begin{definition}
A point $p$ of $X$ is a complex jump point if the dimension of $H_{p}$ is 
$2n+2$.
\end{definition}
(Lai \cite{Lai}
    used the term RC-singular point and Wells \cite{Wells} used the term 
    nongeneric point).

We denote the set of complex jump points by $\cal {J}$.
Then $X-\cal {J}$ has a natural structure of a codimension two CR-manifold.

\begin{definition}
A real $2n+k$ submanifold $X$ in $C^{n+k}$ is an embedded
 CR manifold of codimension $k$ if for all points $p$ in $X$, the 
 complex tangent space $H_{p}$ has real dimension $2n$.
\end{definition}

There is an abstract notion of a CR structure, namely:
\begin{definition}
A real $2n+k$ manifold $X$ will be a codimension $k$ CR manifold if
there is a complex subbundle $L$ of the complexified tangent bundle
$C \otimes TM$ such that
$[L,L] \subset L$ and $L\bigcap \overline{L} = 0$.
\end{definition}
All embedded CR manifolds  are CR manifolds, simply by 
identifying the subbundle $L$ in the latter definition with the $i$ 
eigenbundle 
 $H^{10}$ of the map $J$ for the complexified  bundle $C\otimes H$.
The lion's share of the work on CR structures has been on 
trying to determine when a CR structure can be realized as a real 
submanifold of a complex space.  We will not be concerned here with 
those questions.

\subsection{Nash Blow-ups}

Nash blow-ups are a technique for trying to resolve
 singularities of 
 embedded varieties.   It is unknown whether or 
 not repeated applications of Nash blow-ups will resolve all 
 singularities.   We will look at an example of how to use the Nash blow-up 
to resolve a node of a plane curve. Consider the plane curve $X$ given as 
the zero locus of the polynomial $f(x,y) = y^{2}-x^{3}-x^{2}$.  Since 
both partials are zero at the origin, the origin is a singular 
point.   The Gauss map
$$\sigma:X-(0,0) \rightarrow P^{1},$$
where $P^{1}$ denotes the complex projective line, is defined by 
sending each point of $X-(0,0)$ to its tangent line.  Thus 
$$\sigma(p) = (\frac{\partial f}{\partial y}:-\frac{\partial 
f}{\partial x})= (2y: 2x + 3x^{2}).$$
The Nash blow-up is the closure of this graph in $X\times P^{1}$.  For 
this example, it can be explicitly checked using local coordinates that 
the closure is smooth, with two points sitting over 
the origin $(0,0)$, namely the points $(0,0) \times (1:1)$ and 
 $(0,0) \times (1:-1)$, reflecting that for this plane curve the 
 lines $x=y$ and $x=-y$  are the natural tangents at the origin.

For more information on Nash 
blow-ups, see \cite{Harris}, page 221. It should be noted that the 
Nash blow-up is not the same as the usual blow-up. 

\section{Lai's Work}

The major work on the global properities of embedded CR structures 
has so far been done by Lai in \cite{Lai}.  (See also the work of 
Webster in \cite{Webster1}\cite{Webster2}\cite{Webster3}  and
Coffman in \cite{Coffman1} \cite{Coffman2} \cite{Coffman3}).  
Since we use his work as a 
springboard for this paper, we quickly review his results and 
techniques.  He concentrates on the Gauss map
$$\sigma: X \rightarrow Gr(2n+2, C^{n+2})$$ which maps each point $p 
\in X$ to its tangent space $T_{p}X$ in the Grassmannian $Gr(2n+2, 
C^{n+2})$.
Set 
$${\cal C} = \{ \Lambda \in Gr(2n+2, C^{n+2}): \Lambda \;\mbox{inherits a 
complex structure from}\; C^{n+2}\}.$$
Since generic elements in $Gr(2n+2, C^{n+2})$ will not be themselves 
complex spaces (but instead will only contain a complex subspace of 
real dimension $2n$), $\cal C$ will be a proper subvariety in $Gr(2n+2, 
C^{n+2})$.
The next lemma follows from the definitions:

\begin{lemma}
A point $p$ in $X$ will be a complex jump point precisely when 
$\sigma(p)\in \cal C$.
\end{lemma}
Lai describes the cycle corresponding to 
$\cal C$ in terms of the special Shubert cycles (which generate 
the ring structure of the homology $H_{*}(Gr(2n+2, C^{n+2}))$.  
By pulling back the information from the Grassmannian, Lai 
showed in \cite{Lai}:
\begin{theorem}[Lai]
 Let $F$ be a real k-dimensional manifold and $M$ a real
  2n-dimensional almost complex manifold. Let $i: F \rightarrow  M$
   be an immersion. Assume $2n-2=k$. Then
$$\Omega(F) + \sum_{r=0}^{n-1}\bar{\Omega}(F)^{n-r-1}\cup i^{*}c_{r}(M) =
2\sigma ^{*}(\sigma (F) \cdot {\cal C}).$$
\end{theorem}
Here $\Omega(F)$ is the Euler class of $F$,  $\bar{\Omega}(F)$
  is the Euler class of the normal bundle of $F$ in $M$, $\sigma$ 
  and $\cal C$ are analogues of our earlier definitions, $\cup$ is 
  the cup product and 
  $\sigma ^{*}(\sigma (F) \cdot {\cal C})$ denotes
   the pullback of  $\sigma (F) \cdot {\cal C}$, which is the 
   Poincare dual of the  
   intersection
  product of $\sigma (F)$ and {\cal C} in $H_{*}(Gr(2n+2, C^{n+2}))$. 
     In our case 
  the manifold $M$ is simply $C^{n+2}$ and the submanifold 
  $F$ is $X$.  Note that the right-hand side of this formula is an 
  algebraic count of the number of complex jump points, showing
  that there are topological reasons for the 
  existence of jump points.

  The initial part of Lai's proof needs to use that for generic $X$, 
  the image $\sigma(X)$ will transversally intersect the subvariety $\cal 
  C$. The assumption of transversality will be seen to be the needed 
  condition for the CR-Nash blow-up to be smooth.

  In the case when $k=2n-2$ (the codimension two case),
   we have that $\sigma (X) \cap \cal C$ will be a 
  finite number of points.  Thus in codimension two, there are 
  generically only a finite number of complex jump points.
  
  \section{Flags and the CR-Nash Blow-up}
  For this section, we will denote a complex $n$ dimensional subspace 
  by $\Sigma$ and a real $2n + 2$ dimensional subspace by $\Lambda$.
    Set
  $$F = \{ (\Sigma, \Lambda):\Sigma \subset \Lambda \subset 
  C^{n+2}\}.$$
  F is an example of a flag manifold. By similar argument as in 
  \cite{Harris} in example 11.40, F is locally isomorphic to the 
  product $Gr_{C}(n, n+2) \times Gr(2n, 2n+2)$, where $Gr_{C}(n, 
  n+2)$ is the Grassmannian of complex $n$ dimensional subspaces of 
  the complex space $C^{n+2}$.  Note that there is a natural map from 
  $F$ to $Gr(2n+2, C^{n+2})$,  given by simply sending each $(\Sigma, 
  \Lambda)$ to $\Lambda$.  The inverse image of the map over any 
  $\Lambda \not\in \cal C$ will be a single point, but over a $
\Lambda \in \cal C$ the inverse image will be the full complex
Grassmannian   $Gr_{C}(n, n+1)$.

There are natural universal bundles over a flag, analogous to the 
universal bundles for Grassmannians.  Let $U_{n}$ be the complex rank $n$  
vector bundle whose fibre over a point $(\Sigma, \Lambda)$ consists of 
points in $\Sigma$.  This bundle is a subbundle of the real rank $2n+2$ 
vector bundle $U_{2n+2}$, whose fibre over the point $(\Sigma, 
\Lambda)$ consists of the points in $\Lambda$.

We now want to extend the Gauss map.

\begin{definition}
The CR-Gauss map $\tau: X - J \rightarrow F$ is the map
$$\tau(p) = (H_{p},T_{p}X).$$
\end{definition}
Note that the pullback of the vector bundle $U_{n}$ is the vector 
bundle $H$ and the pullback of the vector bundle $U_{2n+2}$ is the 
tangent bundle $TX$.  Also,
the CR-Gauss map is not defined at complex jump points, 
since at these points, $H_{p}$ is the full tangent space $T_{p}X$.
\begin{definition}
The CR-Nash blow-up $\tilde{X}$ is the closure of the graph of the CR-Gauss map in 
the space $X \times F$.
\end{definition}
This is the CR analogue of the traditional Nash blow-up.  

We can now state the main theorm of this paper.

\begin{theorem}
Let $X$ be a real $2n+2$ dimensional submanifold of the complex 
space  $C^{n+2}$ such that the image of $X$ under the Gauss map 
$\sigma$ intersects transversally the subvariety $\cal C$ in the real 
Grassmannian $Gr(2n+2, C^{n+2})$.  Then the CR-Nash blow-up $\tilde{X}$ is 
a smooth manifold.
\end{theorem}

\section{Transversality in local coordinates}
In order to prove the main theorem we must first have a good 
description of when the image of the Gauss map of $X$ intersects $\cal 
C$ transversally.  As is common with Grassmannians, we will dualize 
the Gauss map, now defining it as:
$$\sigma:X \rightarrow Gr(2,C^{n+2}),$$
with $\sigma(p) = N_{p}$, the conormal bundle. The analogue of the 
subvariety of $2n+2$ planes that inherit a complex structure from 
$C^{n+2}$ will be
$${\cal C} = \{ \Lambda \in  Gr(2,C^{n+2}): \Lambda \; \mbox{inherits a 
complex structure from}\; C^{n+2}\}.$$

Viewing $C^{n+2}$ as the real vector space $R^{2n+4}$, complex 
conjugation becomes a linear map $J:R^{2n+4} \rightarrow R^{2n+4}$ 
with $J^{2} = -I.$ Extending the map $J$ to $C \otimes R^{2n+4}$ 
allows us to split $C \otimes R^{2n+4}$ into its $+i$ and $-i$ 
eigenspaces, which are denoted $H^{10}$ and $H^{01}$ respectively:
$$C \otimes R^{2n+4} = H^{10} \oplus H^{01}.$$  For a vector $v \in C 
\otimes R^{2n+4}$, we write this splitting as 
$$v = v^{10} \oplus v^{01} = (v^{10},v^{01}).$$  Following from the 
discussion in section 3.2 of Boggess\cite{Boggess}, we can show:
\begin{lemma}
Two vectors $v$ and $w$ in $C \otimes R^{2n+4}$  will span a two-plane 
in $\cal C$ if $v\wedge w \neq 0$ but 
$$v^{10}\wedge w^{10} = 
v^{01}\wedge w^{01} = 0.$$
\end{lemma}

We will need to understand $\cal C$'s local coordinates with respect 
to the various  
coordinate systems for the Grassmannian $Gr_{C}(2, C^{2n+4})$ given 
by the Plucker embedding of $Gr_{C}(2, C^{2n+4})$
into the complex projective space $P^{2(2n+2) -1}$.  Recall how this 
map is defined.  Let vector $v$ and $w$ span the two-plane $\Lambda$.  
Then the Plucker embedding is given by $v\wedge w$. If we choose a 
basis for $C\otimes R^{2n+4}$ and use the splitting $H^{10} \oplus 
H^{01}$, we can write each two-plane as the span  of the two 
rows:
$$\pmatrix{v\cr w \cr} = \pmatrix{v^{(10)}&v^{(01)}\cr w^{(10)}&w^{(01)}\cr}=
\pmatrix{v_{1}&\dots &v_{2n+4}& v_{\overline{1}}&\dots &v_{\overline{2n+4}}\cr 
w_{1}&\dots &w_{2n+4}& w_{\overline{1}}&\dots &w_{\overline{2n+4}}\cr}.$$
Then the Plucker embedding is given by the determinants of the two by two 
minors in the above 
matrix.  This is not yet a coordinate system.  At least one of 
these determinants must be nonzero. Here we will assume  that the 
two by two minor  
$$\pmatrix{v_{n+2}&v_{\overline{n+2}}\cr  
w_{n+2}&w_{\overline{n+2}}\cr}$$
is invertible.  By a change of basis of $C^{2n+4}$ we can in fact 
assume that
$$\pmatrix{v_{n+2}&v_{\overline{n+2}}\cr  
w_{n+2}&w_{\overline{n+2}}\cr} = 
\pmatrix{1&1\cr i&-i\cr}.$$
By keeping this matrix fixed and then considering the Plucker 
embedding, we obtain a coordinate system on 
the open set in the  complex Grassmannian where
 $\pmatrix{v_{n+2}&v_{\overline{n+2}}\cr  
w_{n+2}&w_{\overline{n+2}}\cr} = 
\pmatrix{1&1\cr i&-i\cr}.$
  Then the coordinates on this open set 
  for the complex Grassmannian will be given by 
  $$ u_{k,n+2}=i v_{k} - w_{k},$$
  (the $(k,n+2)$ parts of the wedge product), 
  $$u_{k,\overline{n+2}} = -i v_{k} - w_{k},$$
  (the $(k,\overline{n+2})$ parts of the wedge product),
  $$u_{\overline{k},n+2}=i v_{\overline{k}} - w_{\overline{k}},$$ 
  (the $(\overline{k},n+2)$ parts of the wedge product) and
   $$u_{\overline{k},\overline{n+2}}=-i v_{\overline{k}} - 
   w_{\overline{k}},$$
    (the 
   $(\overline{k},\overline{n+2})$ parts of the wedge product).
   
   \noindent  On our fixed open subset of the 
Grassmannian, $\cal C$ will be the linear subvariety
$$u_{k,n+2} = u_{\overline{k},\overline{n+2}} = 0,$$
since  $\cal C$ is where $v^{10}\wedge w^{10} = 
v^{01}\wedge w^{01} = 0.$ (Note that this shows that the dimension of 
$\cal C$ is $2n+2$).
  Fix the  basis for the tangent space to 
the whole Grassmannian, on our open subset, to be
 $\frac{\partial}{\partial u_{k,n+2}},
\frac{\partial}{\partial u_{k,\overline{n+2}}},
\frac{\partial}{\partial u_{\overline{k},n+2}},
\frac{\partial}{\partial u_{\overline{k},\overline{n+2}}}$.
The tangent space to $\cal C$ is the span of the vectors
$\frac{\partial}{\partial u_{k,\overline{n+2}}},
\frac{\partial}{\partial u_{\overline{k},n+2}}$.  Then in terms of this basis,
we can describe the tangent space to $\cal C$ by the 
the $(2n+2) \times (4n+4)$ matrix
$$\pmatrix{0&I&0&0\cr  
0&0&I&0\cr} ,$$
where each $I$ is an $(n+1) \times (n+1)$ identity matrix.  Here the 
first $(n+1)$ columns correspond to the $(k,n+2)$ parts of the wedge 
product, the next $(n+1)$ columns
 correspond to the $(k,\overline{n+2})$ parts of the wedge product, 
 etc.  The first $n+1$ rows correspond to $\cal C$'s tangent vectors
 $\frac{\partial}{\partial u_{k,\overline{n+2}}}$ and the last $n+1$ 
 rows correspond to $\cal C$'s tangent vectors
  $\frac{\partial}{\partial u_{\overline{k},n+2}}$.

Return to our manifold $X$.  At a point $p\in X$, we can describe $X$ 
as the zero locus of two smooth real-valued functions:
$$X = (\rho_{1}=0)\cap (\rho_{2}=0).$$
Then the Gauss map will be:
$$\sigma(x) = \mbox{span}(d\rho_{1},d\rho_{2}).$$
Then a complex jump point (those points whose image under $\sigma$ 
lands in $\cal C$) will be those points where 
$\partial\rho_{1}\wedge\partial\rho_{2} = 0$ (see section 7.1, lemma 4 
in \cite{Boggess}).

We want to find clean conditions for when the intersection of 
$\sigma(X)$ with $\cal C$ is transverse.  Thus we must look at the 
Jabobian $D\sigma$.  Let $p\in X$ be a complex jump point.  Change 
coordinates so that $p$ is the origin in $C^{n+2}$. Rotate the coordinate 
system so that locally, about the origin, $X$ is the zero locus of the 
two smooth functions
$$\rho_{1} = z_{n+2} + \overline{z_{n+2}} + f_{1}$$
$$\rho_{2} = i(z_{n+2} - \overline{z_{n+2}}) + f_{2},$$
where the functions $f_{1}$ and $f_{2}$ are smooth functions that 
vanish to second order 
at the origin. 
  Since we have $d\rho_{1}(0) = dz_{n+2} + d\overline{z_{n+2}}$ and 
  $d\rho_{2}(0) = i(dz_{n+2} - d\overline{z_{n+2}}),$ the origin does 
  map to a point in $\cal C$.  Both $X$ and $\cal C$  have real 
   dimension $2n + 2$, which is half of the dimension of the ambient 
   Grassmanian.  Thus we will have a transverse intersection if the 
   respective tangent spaces span the full tangent space of the 
   Grassmanian.

The Plucker coordinates of the Gauss map for $X$ are given by the two 
by two minors of the matrix 
$$\pmatrix{\partial\rho_{1}&\overline{\partial}\rho_{1}\cr 
\partial\rho_{2}& \overline{\partial}\rho_{2}\cr}$$
and hence are
$$ u_{k,n+2} = i\frac{\partial \rho_{1}}{\partial z_{k}} -
\frac{\partial \rho_{2}}{\partial z_{k}}, $$
$$ u_{k,\overline{n+2}}= -i \frac{\partial \rho_{1}}{\partial z_{k}}-
\frac{\partial \rho_{2}}{\partial \overline{z}_{k}}, $$ 
$$u_{\overline{k}, n+2} = 
i \frac{\partial \rho_{1}}{\partial \overline{z}_{k}}
 - \frac{\partial \rho_{2}}{\partial z_{k}},$$
 and
 $$u_{\overline{k},\overline{n+2}}=
  -i\frac{\partial \rho_{1}}{\partial \overline{z}_{k}} -
\frac{\partial \rho_{2}}{\partial \overline{z}_{k}}.$$
To compute the Jacobian, we need to differentiate this map with 
respect to a local coordinate system of $X$.  
  We can assume that at the origin  the local 
coordinate system for $X$ is given by 
$z_{1},\ldots,z_{n+1},\overline{z_{1}},\ldots,\overline{z_{n+1}}$.
 Then the tangent space 
to the image at $X$ will be the $(2n+2)\times (4n+4)$ matrix

$$\pmatrix{\frac{\partial}{\partial z_{1}}
(i\frac{\partial \rho_{1}}{\partial z_{k}} -
\frac{\partial \rho_{2}}{\partial z_{k}})&
\frac{\partial}{\partial z_{1}}(-i \frac{\partial \rho_{1}}{\partial z_{k}}-
\frac{\partial \rho_{2}}{\partial \overline{z}_{k}})
 &\frac{\partial}{\partial z_{1}}
(i \frac{\partial \rho_{1}}{\partial \overline{z}_{k}}
 - \frac{\partial \rho_{2}}{\partial z_{k}})
&\frac{\partial}{\partial z_{1}}
(-i\frac{\partial \rho_{1}}{\partial \overline{z}_{k}} -
\frac{\partial \rho_{2}}{\partial \overline{z}_{k}}) \cr
\vdots&\vdots&\vdots&\vdots
 \cr  \frac{\partial}{\partial \overline{z}_{n+1}}
(i\frac{\partial \rho_{1}}{\partial z_{k}} -
\frac{\partial \rho_{2}}{\partial z_{k}})&
\frac{\partial}{\partial \overline{z}_{n+1}}
(-i \frac{\partial \rho_{1}}{\partial z_{k}}-
\frac{\partial \rho_{2}}{\partial \overline{z}_{k}})
 &\frac{\partial}{\partial \overline{z}_{n+1}}
( i \frac{\partial \rho_{1}}{\partial \overline{z}_{k}}
 - \frac{\partial \rho_{2}}{\partial z_{k}})
&\frac{\partial}{\partial \overline{z}_{n+1}}
(-i\frac{\partial \rho_{1}}{\partial \overline{z}_{k}} -
\frac{\partial \rho_{2}}{\partial \overline{z}_{k}})\cr}.$$
Here the $k$ are running from $1$ to $n+1$.
Using our earlier description of the tangent space of $\cal C$, we see 
that transversality will occur when the  
$(2n+2) \times (2n+2)$ minor of the above 
matrix formed from the first $(n+1)$ columns 
and the 
last $(n+1)$ columns of the above matrix is invertible.

\section{Smoothness}

We now want to prove the main theorm of this paper, namely

\noindent{\it Let $X$ be a real $2n+2$ dimensional submanifold of the complex 
space  $C^{n+2}$ such that the image of $X$ under the Gauss map 
$\sigma$ intersects transversally the subvariety $\cal C$ in the real 
Grassmannian $Gr(2n+2, C^{n+2})$.  Then the CR-Nash blow-up $\tilde{X}$ is 
a smooth manifold.}

\noindent We will reduce this to the   standard blow-up of the origin in 
$C^{n+1}$ (as in \cite{Griffiths-Harris}, page 182), which is well 
known to be smooth.

In a manner similar to example 11.40 in \cite{Harris}, we can locally 
write our flag manifold $F$ as sitting inside $Gr(2n, 2n + 2 )\times 
Gr(2n+2, 2n+4)$.  The CR-Gauss map $\tau$  projected onto the second factor is 
the traditional Gauss map.  Since our manifold $X$ is smooth in 
$C^{n+2}$, this part of the closure of $\tau(X-\cal {J})$ will be smooth.  
The part where the closure can fail to be smooth will be the part of 
$\tau$ that is projected onto the first factor.  Since
 $\tau(p)= (H_{p},T_{p}X)$ at non-jump points $p$, it is 
the first factor $H_{p}$ that fails to be defined at jump points and 
is the source of the difficulties.

Let $p$ be an isolated jump point at which  the Gauss map 
$\sigma$ intersects transversally the subvariety $\cal C$.  We know 
that at this point the tangent space $T_{p}X$ inherits a complex 
structure from the ambient space and can thus be identified to $C^{n+1}$.
Then our flag can be identified with 
$Gr_{C}¥(n, C^{n+1} )\times 
Gr(2n+2, 2n+4)$, where $Gr_{C}¥(n, C^{n+1})$ is the Grassmannain of 
complex subspace of dimension $n$ in $C^{n+1}$.  At points $q$ near 
$p$, we know that $\partial \rho_{1}(q) \wedge
 \partial \rho_{2}(q)\neq 0$ (which of course via duality defines the 
 subspace $H_{q}$) but   
  $\partial \rho_{1}(p) \wedge
 \partial \rho_{2}(p)= 0$.

  Using the notation from the previous section, we know that the 
  Plucker coordinates of the Gauss map of $X$ are:

$$ u_{k,n+2} = i\frac{\partial \rho_{1}}{\partial z_{k}} -
\frac{\partial \rho_{2}}{\partial z_{k}} $$
$$ = i\frac{\partial f^{1}}{\partial z_{k}} -
\frac{\partial f^{2}}{\partial z_{k}} $$
and

$$u_{\overline{k},\overline{n+2}}=
  -i\frac{\partial \rho_{1}}{\partial \overline{z}_{k}} -
\frac{\partial \rho_{2}}{\partial \overline{z}_{k}}$$
$$=-i\frac{\partial f^{1}}{\partial \overline{z}_{k}} -
\frac{\partial f^{2}}{\partial \overline{z}_{k}}.$$
By the transversality assumption, we have that the $(2n+2)\times 
(2n+2)$ matrix,

$$\pmatrix{\frac{\partial}{\partial z_{1}}
(i\frac{\partial \rho_{1}}{\partial z_{k}} -
\frac{\partial \rho_{2}}{\partial z_{k}})
&\frac{\partial}{\partial z_{1}}
(-i\frac{\partial \rho_{1}}{\overline{z}_{k}} -
\frac{\partial \rho_{2}}{\overline{z}_{k}}) \cr
\vdots&\vdots
 \cr  \frac{\partial}{\partial \overline{z}_{n+1}}
(i\frac{\partial \rho_{1}}{\partial z_{k}} -
\frac{\partial \rho_{2}}{\partial z_{k}})
&\frac{\partial}{\partial \overline{z}_{n+1}}
(-i\frac{\partial \rho_{1}}{\overline{z}_{k}} -
\frac{\partial \rho_{2}}{\overline{z}_{k}})\cr}
= \pmatrix{\frac{\partial u_{k,n+2}}{\partial z_{1}}&
\frac{\partial u_{\overline{k},\overline{n+2}}}{\partial z_{1}}\cr
\vdots&\vdots
 \cr \frac{\partial u_{k,n+2}}
 {{\partial \overline{z}_{n+1}}}&
\frac{\partial u_{k,n+2}}{{\partial \overline{z}_{n+1}}}\cr},$$
where $k= 1, \ldots, n+1$,is invertible.  
Then we can choose a (real) coordinate system 
 $w_{1},\ldots
 , w_{2n+2}$ for $X$ so that
 
 $$u_{k,n+2} = w_{k}+ iw_{n+k} + \mbox{higher order terms}$$
and 
$$u_{\overline{k},\overline{n+2}}
 = w_{k}-i w_{n+k}+ \mbox{higher order terms}.$$
 
 Let $\bigwedge^{(2,0)}C^{n+2}$ denote the vector space of 
 $(2,0)$ forms on $C^{n+2}$.  
 There is the natural map
 $$X \rightarrow \mbox{$\bigwedge^{(2,0)}$}C^{n+2}$$
 given by sending a point $q$ to $\partial \rho_{1}(q) \wedge
 \partial \rho_{2}(q)$.  Away from the complex jump points, we have 
 the map
 
 $$X - J \rightarrow P(\mbox{$\bigwedge^{(2,0)}$}C^{n+2}),$$
 where $P(\bigwedge^{(2,0)}C^{n+2})$ denotes the projectivization of 
 $\bigwedge^{(2,0)}C^{n+2}$.  We want to look at the closure of this 
 graph in $X \times P(\bigwedge^{(2,0)}C^{n+2})$.  By our choice of 
 local coordinates, we have
 $$\partial \rho_{1} \wedge
 \partial \rho_{2} = \sum (w_{k}+iw_{n+k})dz_{k}\wedge dz_{n+2} + \;
 \mbox{higher order terms}$$
 and 
 $$\overline{\partial} \rho_{1} \wedge
 \overline{\partial} \rho_{2} = \sum (w_{k}-iw_{n+2+k})d\overline{z}_{k}
 \wedge d\overline{z}_{n+2} + \;
 \mbox{higher order terms}.$$
 
 But then the closure will be smooth, since up to higher order we can 
 view the map as a map 
 $X \rightarrow P^{n}$
 given by 
 $$(w_{1}+iw_{1+n+2},\ldots, w_{n+1}+iw_{n+1 + n+2})\rightarrow
 (w_{1}+iw_{1+n+2}:\ldots: w_{n+1}+iw_{n+1 + n+2})$$
 and thus the closure is smooth (again, this is known and can also be 
 directly calculated).
  Under duality, we have that the graph in $X \times G_{C}(n, n+1)$ 
  will be smooth, completing the proof.

\section{Extending the Levi form to the blow-up}

The key tool for understanding CR structures is the Levi form, which 
is a vector-valued map:
$$L: = H^{10} \times H^{01}: \rightarrow C\otimes TX/(H^{10} \oplus  H^{01}),$$
defined as follows.  Let $p\in X$ and let $v_{p}\in H^{10}_{p}$ and 
$w_{p}\in H^{01}_{p}$.  Extend $v_{p}$ to a vector field $v$ in 
$H^{10}$ and $w_{p}$ to a vector field $w$ in $H^{01}$.  Then define 
$L(v_{p}, w_{p})$ as
$L(v_{p}, w_{p}) = \pi_{p}[v,w],$
where $[v,w]$ is the Lie bracket and 
$\pi_{p}: C\otimes TX \rightarrow C\otimes TX/(H^{10} \oplus  H^{01})$
 is the natural 
projection map.  Here we are using that the Lie bracket of two 
tangent vectors is again a tangent vector and that there is a natural 
projection map to $C \otimes TX/(H^{10} \oplus  H^{01})$.  At complex jump points, 
the Levi form will be undefined due to the lack of the natural 
projection map.  

There is an alternative approach for defining the Levi form. Again, we 
restrict attention to where $X$ has a CR structure.  As before, $X$ 
is locally defined in ${C}^{n+2}$ as the zero locus  of the 
 functions $\rho_{1}$ and $\rho_{2}$ but now assume that the 
 vectors $\bigtriangledown \rho_{1}$ and $\bigtriangledown \rho_{2}$ 
 form an orthonormal basis for the normal bundle $N$.  (We will be 
 using throughout the natural Hermitian metric on $C^{n+2}$, allowing 
 us to identify various bundles and their dual spaces, an 
 identification that will usually not be explicitly made).  Using that the 
 normal bundle $N$ is isomorphic to the
  bundle $C \otimes TX/(H^{10} \oplus  H^{01})$, under the map $J$, we 
 can define the Levi form as follows.  Let 
 $$v = \sum_{j=1}^{n+2} v_{j} \frac{\partial}{\partial z_{j}}$$
 be vector in  $H^{10}$ and 
 $$w = \sum_{j=1}^{n+2} v_{\overline{j}} \frac{\partial}{\partial 
 \overline{z}_{j}}$$
 a vector in $H^{01}$.  Then the map
 $$\tilde{L}: = H^{10} \times H^{01}: \rightarrow C\otimes N$$
 defined by 
$$\tilde{L}(v,w) = -(\sum_{j,k=1}^{n+2} 
\frac{\partial^{2}\rho_{1}}{\partial z_{j} \partial 
\overline{z}_{k}}v_{j}w_{\overline{k}})\bigtriangledown \rho_{1} + 
\sum_{j,k=1}^{n+2} 
\frac{\partial^{2}\rho_{2}}{\partial z_{j} \partial 
\overline{z}_{k}}v_{j}w_{\overline{k}})\bigtriangledown \rho_{2}),$$
is equivalent to the Levi form, as shown in \cite{Boggess} in section 
10.2.

We want to extend this to the CR-Nash blow-up $\tilde{X}$. 
A point 
in $\tilde{X}$ is described by specifying a point $p\in X$ and a $2n$ 
dimensional subspace $H^{10} \oplus  H^{01}$ of 
$C \otimes TX$ and thus as $(p, H^{10} \oplus  H^{01}, C\otimes TX)$ 
in $X\times F$.  Over the flag $F$ there are the natural universal 
bundles $C\otimes U_{n}$ and $C\otimes U_{2n+2}$ 
which match up, away from the complex jump points of $X$, 
with the bundles $H^{10} \oplus  H^{01}$ and
$C \otimes TX$, respectively.   Further the isomorphism from 
the normal bundle $N$ 
(which is ${C}^{n+2}/ 
{C}\otimes TX$)
 to ${C}\otimes TX/(H^{10} \oplus  H^{01})$  extends, which we 
 will still denote by $J$.
 \begin{definition}
Let $(p, H^{10} \oplus  H^{01}, C\otimes TX)$ be a point in the CR-Nash 
blow-up of $X$.  Let $v_{p} = 
\sum_{j=1}^{n+2} v_{j} \frac{\partial}{\partial z_{j}}\in H^{10}_{p}$ and 
$w_{p} = \sum_{j=1}^{n+2} v_{\overline{j}} \frac{\partial}{\partial 
 \overline{z}_{j}}\in H^{01}_{p}$.    Define the 
Levi form to be the map
$$L: H^{10} \times H^{01}: \rightarrow C\otimes 
TX/(H^{10} \oplus  H^{01})$$
given by
$$L(v_{p}, w_{p}) = J(-(\sum_{j,k=1}^{n+2} 
\frac{\partial^{2}\rho_{1}}{\partial z_{j} \partial 
\overline{z}_{k}}v_{j}w_{\overline{k}})\bigtriangledown \rho_{1} + 
\sum_{j,k=1}^{n+2} 
\frac{\partial^{2}\rho_{2}}{\partial z_{j} \partial 
\overline{z}_{k}}v_{j}w_{\overline{k}})\bigtriangledown \rho_{2})).$$
\end{definition}

% This problem does not arise for the CR-Nash blow-up $\tilde{X}$. A point 
% in $\tilde{X}$ is described by specifying a point $p\in X$ and a $2n$ 
% dimensional subspace $H^{10} \oplus  H^{01}$ of 
% $C \otimes TX$ and thus as $(p, H^{10} \oplus  H^{01}, C\otimes TX)$ 
% in $X\times F$.  Over the flag $F$ there are the natural universal 
% bundles which match up, away from the complex jump points of $X$, 
% with the bundles $H^{10} \oplus  H^{01}$ and
% $C \otimes TX$. Thus we have the following exact sequence of vector 
% bundles over the blow-up $\tilde{X}$:
% $$ 0 \rightarrow H^{10} \oplus  H^{01} \rightarrow C\otimes TX 
% \rightarrow C\otimes TX/(H^{10} \oplus  H^{01} \rightarrow 0.$$
% \begin{definition}
% Let $(p, H^{10} \oplus  H^{01}, C\otimes TX)$ be a point in the CR-Nash 
% blow-up of $X$.  Let $v_{p}\in H^{10}_{p}$ and 
% $w_{p}\in H^{01}_{p}$.  As tangent vectors on $X$,
%  extend $v_{p}$ to a vector field $v$ in 
% $H^{10}$ and $w_{p}$ to a vector field $w$ in $H^{01}$.  Define the 
% Levi form to be the map
% $$L: H^{10} \times H^{01}: \rightarrow C\otimes 
% TX/(H^{10} \oplus  H^{01})$$
% given by
% $$L(v_{p}, w_{p}) = \pi_{p}[v,w],$$
% where $[v,w]$ is the Lie bracket of tangent vector on $TX$ and 
% $\pi_{p}: C\otimes TX \rightarrow C\otimes TX/(H^{10} \oplus  H^{01})$
%  is the natural 
% projection map.
% \end{definition}
% 
% This can be shown by the usual arguments that the Levi form is 
% independent fo the extensions $v$ and $w$. 

\section{An example of a global obstruction: Levi non-degeneracy}

The Levi form has been the main tool in trying to solve the local 
equivalence problem for CR structures;
much of the previous work depended on placing various algebraic restrictions 
on the Levi form. We will find topological obstructions for  the 
Levi form to be nondegenerate.  The same obstructions will be seen to 
effect the local work in \cite{Mizner}.  

Locally on the Nash blow-up $\tilde{X}$, choose sections for $ H^{10}$ 
(which will give us sections for $H^{01}$), 
 and $TX/(H^{10} \oplus  H^{01})$.  Then the Levi form becomes 
two $n \times n$ Hermitian matrices $(L_{1},L_{2})$.  Consider the 
degree $n$ homogeneous
polynomial (first introduced by Mizner \cite{Mizner}):
$$P(x,y) = det(x L_{1} + y L_{2}).$$
If we change the choice of sections for $ H^{10}$ by an element 
$g\in GL(n, C)$, the polynomial is altered by multiplying all of its 
coefficients by the factor $|det(g)|^{-2}$. 
 Changing sections for $TX/(H^{10} \oplus  H^{01})$ will correspond to 
 making a homogeneous change of coordinates of the polynomial $P(x,y)$.
 Thus the polynomial $P(x,y)$ can be viewed as a section of the 
 bundle $\wedge^{n}H^{01 *} \otimes \wedge^{n}H^{01 *}\otimes 
 S^{n}TX/(H^{10} \oplus  H^{01})$, 
 where $S^{n}$ denotes the $(n)th$ symmetric product 
 of $TX/(H^{10} \oplus  H^{01})$.
 
 We will concentrate on determining 
 the topological obstructions that would force the polynomial $P(x,y)$
 to be the zero polynomial (which means that the two Hermitian 
 matrices $L_{1}$ and $L_{2}$ share a nontrivial element in their 
 kernels).  From 20.10.5 in \cite{Bott-Tu}, we see that a complex 
 vector bundle has a non-vanishing section when its top Chern class is 
 zero.   Since $\wedge^{n}H^{01 *} \otimes \wedge^{n}H^{01 *}\otimes 
 S^{n}W$ has rank $n+1$, if 
 $$c_{n+1}(\wedge^{n}H^{01 *} \otimes \wedge^{n}H^{01 *}\otimes 
 S^{n}W) \neq 0,$$
 then there must be points on the Nash blow-up at which the 
 polynonial $P(x,y)$ is the zero polynomial. 
 
 Now to see how the vanishing of the polyomial $P(x,y)$ relates to 
 Levi non-degeneracy.
 \begin{definition}
 A Levi form $L=(L_{1},L_{2})$ is non-degenerate if
 
 i.  $L_{1}$ and $L_{2}$ are linearly independent.
 
 ii.  $L_{1}$ and $L_{2}$ do not share a common nonzero kernel.
  \end{definition} 
 This has been an important idea  in the work of many of the people 
 mentioned in the introduction.  Note that
  if $L_{1}$ and $L_{2}$ do  share a common nonzero kernel, then 
 $P(x,y)$ is the zero polynomial.  Thus if 
 $c_{n+1}(\wedge^{n}H^{01 *} \otimes \wedge^{n}H^{01 *}\otimes 
 S^{n}W) \neq 0$, the Levi form on the blow-up cannot be Levi 
 non-degenerate at every point.

 \section{Questions}
 There should be nothing particularly special about codimension two 
 manifolds.  One can easily define a CR-Nash blow up for any 
 codimensional submanifold of a complex space.  We suspect that if 
 the Gauss map of a submanifold $X$ transversally
intersects the analogue of ${\cal C}$, then the CR-Nash blow up will 
be smooth, for all codimension.  

More difficult is determining if there is a type of CR-Nash blow up 
for an abstract manifold $X$ on which there is a CR structure at most 
points.  If such a blow-up exists, then there is the possibility that 
this will provide 
topological obstructions for embeddibility of 
compact manifolds into a complex space. 

Finally, there is the question of how the work of Harris   
\cite{HarrisG1}, \cite{HarrisG2} and \cite{HarrisG3} on the function 
theory near jump points relates to blow-ups.

\end{document}